\documentclass{vki}

\usepackage{lmodern,pdftexcmds,amsmath}

\usepackage{amsfonts,amsthm,bm} 

\DeclareMathAlphabet{\mathsfbi}{OT1}{\sfdefault}{bx}{sl}
\DeclareMathVersion{sfletters}
\SetSymbolFont{letters}{sfletters}{OML}{ntxsfmi}{b}{it}

\makeatletter
\newcommand{\mathbfsbilow}[1]{%
	\text{\mathversion{sfletters}$\m@th#1$}%
}

\DeclareRobustCommand{\tensor}[1]{%
	\begingroup
	\ifcat\noexpand #1\relax
	\edef\greek@test{\detokenize{#1}}%
	\edef\greek@test{\expandafter\@cdr\greek@test\@nil}%
	\edef\greek@test{\expandafter\@car\greek@test\@nil}%
	\edef\x{\the\lccode\expandafter`\greek@test}%
	\edef\y{\number\expandafter`\greek@test}%
	\ifnum\x=\y\relax
	\mathbfsbilow{#1}%
	\else
	\mathsfbi{#1}%
	\fi
	\else
	\mathsfbi{#1}%
	\fi
	\endgroup
}
\makeatother

\usepackage[most]{tcolorbox}
\usepackage{booktabs}
\usepackage{amsmath}
\usepackage{wrapfig}
\usepackage[linesnumbered,ruled,vlined]{algorithm2e}

\usepackage{mathtools}
\tcbset{enhanced,colback=cyan!2!white,colframe=cyan!90!white,coltitle=blue!20!black,fonttitle=\bfseries}

\usepackage{listings}
\usepackage{hyperref}
\hypersetup{
	colorlinks=true,
	linkcolor=blue,
	filecolor=magenta,      
	citecolor=blue,
	urlcolor=cyan,
}

\begin{document}


\pagenumbering{roman}
\title{Continuous and Discrete LTI Systems}
\author{Miguel Alfonso Mendez\\von Karman Institute for Fluid Dynamics 
\thanks{mendez@vki.ac.be}}

\date{12 February 2020}
\maketitle

This chapter reviews the fundamentals of continuous and discrete Linear Time-Invariant (LTI) systems with Single Input-Single Output (SISO). We start from the general notions of signals and systems, the signal representation problem and the related orthogonal bases in discrete and continuous forms. We then move to the key properties of LTI systems and discuss their eigenfunctions, the input-output relations in the time and frequency domains, the conformal mapping linking the continuous and the discrete formulations, and the modeling via differential and difference equations. Finally, we close with two important applications: (linear) models for time series analysis and forecasting and (linear) digital filters for multi-resolution analysis. This chapter contains seven exercises, the solution of which is provided in the book's webpage\footnote{\url{https://www.datadrivenfluidmechanics.com/download/book/chapter4.zip}}.

\pagenumbering{arabic}
\setcounter{page}{1}
\clearpage{\pagestyle{empty}} 

\tableofcontents
\clearpage{\pagestyle{empty}} 


\section{Signals, Systems, Data and Modeling Fluid Flows}\label{sec4_1}

A \emph{signal} is any function that conveys information about a specific variable (e.g., velocity or pressure) of interest in our analysis. A signal can be a function of one or multiple variables, it can be continuous or discrete, and can have infinite or finite duration/ extension (that is be non-null only within a range of its domain). Signals are produced by \emph{systems}, usually in response to other signals or due to the interaction between different interconnected \emph{subsystems}. In the most general setting, a system is any entity that is capable of manipulating an \emph{input} signal and produce an \emph{output} signal. 

At such a level of abstraction, a vast range of problems in applied science falls within the framework of this chapter. For a fluid dynamicist, the flow past an airfoil is a system in which the inputs are the flow parameters (e.g., free stream-velocity and turbulence) and control parameters (e.g., the angle of attack), and the outputs are the drag and lift components of the aerodynamic force exchanged with the flow. 

Any measurement chain is a system in which the input is the quantity \emph{to be measured}, and the output is the quantity that \emph{is measured}. A hot-wire anemometer, for example, is a complex system that measures the velocity of flow by measuring the heat loss by a wire which is heated by an electrical current. Any signal processing technique for denoising, smoothing and filtering is a system that takes in input the signal to be treated and outputs the signal with modified properties (e.g., with enhanced details or reduced noise).

Regardless of the number of subsystems composing a system, and whether the system is a physical system, a digital replica of it, or an algorithm in a computer program, the relations between input and output are governed by a \emph{mathematical model}. The derivation of such models is instrumental for applications encompassing simulation, prediction or control.
Models can be phrased with various degrees of sophistication, depending on the purposes for which these are developed. They usually take the form of Partial Differential Equations (PDEs), Ordinary Differential Equations (ODEs) or Difference Equations (DEs) or simple algebraic relation. Different models might have different ranges of \emph{validity} (hence different degrees of generalization) and might involve different levels of complexity in their \emph{validation}.

Models can be derived from two different routes, hinging on \emph{data} and \emph{experimentation} in various ways. 
The first route is that of \emph{fundamental principles}, based on the division of a system into subsystems for which empirical observations have allowed to derive well established and validated relations. In the example of the hot-wire anemometer, the subsystems operate according to Newton's cooling law for forced convection, the resistive heating governed by Joule's law, the thermoelectric laws relating the wire resistance to its temperature, and the Ohm's and Kirchhoff's laws governing the electric circuits that are designed to indirectly measure the heat losses. Arguably, none of these laws were formulated with a hot-wire anemometer in mind. Yet, their range of validity is wide enough to accommodate \emph{also} such an application: these laws generalize well. 

In the example of the flow past an airfoil, the system is governed by Navier-Stokes equations, which incorporates other laws such as constitutive relations for the shear stresses (e.g., Newton's law for a Newtonian fluid) and the heat fluxes (e.g., Fourier's Law for conduction). These closure relations led to the notion of fluid properties such as dynamic viscosity and thermal conductivity and were also derived in simple experiments that did not target any specific application.

Because of their remarkable level of generalization, we tend to see these laws as simple mathematical representation of the `laws of nature'. Whether nature is susceptible to mathematical treatment is a question with deep philosophical aspects. As engineers, we accept the pragmatic view of relying on models and laws if these are \emph{validated} and useful. This is the foundation of our scientific and technical background. 

The second route is that of \emph{system identification} or \emph{data-driven approach}, based on the inference of a suitable model from a (usually large) set of input-output data. The model might be constrained to a certain parametric form (e.g., with a given order in the differential equations) or can be completely inferred from data. In the first case, we face a regression problem of identifying the parameters such that the model \emph{fits} the data. In the second case, an algorithm proposes possible models (e.g., in Genetic Programming, see Chapter 14) or uses such a complex parametrization that an analytic form is not particularly interesting (e.g., in Artificial Neural Networks, see Chapter 3). 

This paradigm is certainly not new. An excellent example of system identification is the method proposed by the Swedish Physicist Angstr\"{o}m to measure thermal conductivity of a material in 1861 \citep{Sundqvist1991}. This method consists in using a long metal rod in which a harmonic heat wave is produced by periodically varying the temperature at one end. The spatial attenuation of the heat wave is modeled by the unsteady one-dimensional heat equation, leaving the thermal conductivity as an unknown. A simple formula can be derived for the harmonic response\footnote{The harmonic response of a system is described in Section \ref{sec4_6}.} of the system, and the thermal conductivity is then identified by fitting the model to data. 

Although the data-driven paradigm has a long history \citep{Ljung2008}, its capabilities have been significantly augmented and empowered by our increasing ability to generate enormous amounts of data and by the powerful tools popularized by the ongoing machine learning revolution. Some of these are described in Chapters 1 and 3. The big challenge (and the big opportunity) in big data is to combine physical modeling with the data-driven modeling. It is the author's opinion that the combination of machine learning tools and the general framework of signals and systems can accommodate the formulation of many problems of applied science and, with the required caution, problems in fluid mechanics.

Caution is certainly required as the Navier-Stokes equations governing fluid flows feature the entire spectra of complexities and challenges that a system analyst could think of.
In most configuration of interest, fluid flows are \emph{nonlinear}: a linear combination of input does not produce a predictable linear combination of outputs. It is thus not possible to predict the response of the system from a dictionary of known outputs. Fluid flows are often \emph{high dimensional} and involve many scales: the amount of information required to identify the state of a flow system is tremendous. Fluid flows are often chaotic systems and are thus \emph{unpredictable}: an infinitesimal change in initial or boundary conditions quickly yields very different instantaneous states of the system. Decades of fluid mechanics research aimed at models that are of \emph{reduced order}, which include engineering \emph{closure laws} (e.g., turbulence modeling) or treats the flows in terms of statistical quantities that are usually predictable.

While the Navier Stokes equations are on the top of the ladder of complexity in the signal-system framework, this chapter treats systems that are at the bottom of the ladder: systems that are \emph{linear}, \emph{single dimensional}, \emph{deterministic} and \emph{time invariant}.
These are a subclass of Linear Time (or Translation) Invariant systems (LTI). The interest in such a review is justified by three reasons. The first is that LTI systems are still often encountered in practice, both in physical systems\footnote{see the previous example of Angstr\"{o}m's system.} and in most signal processing operations. The second is that many nonlinear systems can be treated reasonably well as linear systems if the range of operations is close enough fixed state; this is what lead to the remarkable success of many linear control methods described in Chapter 10. The third reason is that this chapter provides the background material to understand time-frequency analysis in Chapter 5 and the multi-resolution analysis underpinning the data decompositions in Chapter 8, the state-space models and the linear control theory in Chapter 10, and the system identification tools in Chapter 12.

Most of the material presented in this chapter can be found in classic textbooks on signal and systems \cite{Oppenheim1996a,Ljung1994,Hsu2013}, signal processing \cite{Williamson1999,Ingle2011,Hayes2011}, orthogonal transforms \cite{Wang2009}, system identification \cite{Oppenheim2015,ljung_book} or control theory \cite{SigurdSkogestad2005,Ogata2009}. We begin this chapter by introducing the relevant notation.

\section{A note about notation and style}\label{sec_4_0}

\begin{figure}[htbp]
	\centering
	\includegraphics[keepaspectratio=true,width=0.9\columnwidth]
	{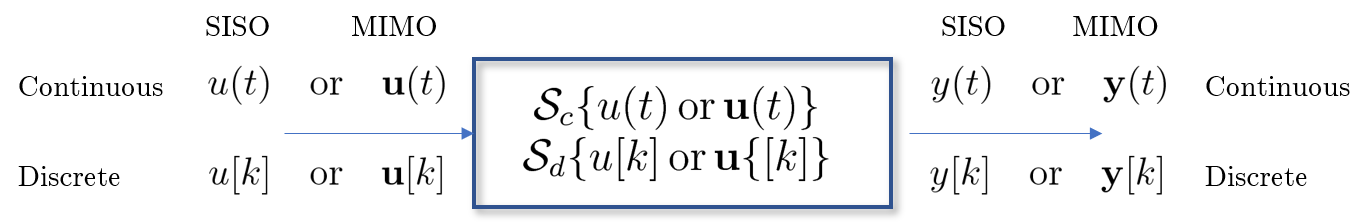}
	\caption{Block diagram rapresentation of systems, classified as Single-Input-Single-Output (SISO) or Multiple-Input-Multiple-Output (MIMO) and continuous or discrete.}\label{fig1}
\end{figure}

In their most general form, continuous and discrete LTI systems admit Multiple Inputs and respond with Multiple Outputs (MIMO systems) or have a Single Input and respond with a Single Output (SISO systems)\footnote{Or can be a combination of the two as in MISO/SIMO systems.}. This classification and the relevant notation is further illustrated in the block diagram in Figure \ref{fig1}.

In a continuous SISO system, input and outputs are denoted respectively as continuous functions $u(t),y(t)\in\mathbb{R}$ with $t\in \mathbb{R}$.
Following the signal processing literature, for discrete systems these are denoted using an index notation, as $u[k], y[k]\in \mathbb{R}$ with $k \in \mathbb{Z}$. In this chapter, discrete and continuous signals are assumed to be linked by a sampling process, hence the time domain is discretized as $t \rightarrow t_k=k\Delta t$ with an index $k\in \mathbb{Z}$, sampling period $\Delta t=1/f_s$, and (constant) sampling frequency $f_s$. Therefore, the notation $u(t_k)$ is equivalent to the index notation $u[k]$, but the second makes no link to the sampling process nor the time axis. In a MIMO system, both the inputs and the outputs are vectors $\mathbf{u}(t),\mathbf{u}[k] \in \mathbb{R}^{n_I}$ and  $\mathbf{y}(t),\mathbf{y}[k] \in \mathbb{R}^{n_O}$, with $n_I$ and $n_O$ the number of inputs and outputs. MIMO systems are better treated in state-space representation presented in Chapter 10, hence this chapter only focuses on SISO systems.

Signals and systems can contain a deterministic and a stochastic part. While the focus is mostly on the deterministic part, the treatment of stochastic signals is briefly recalled in section \ref{sec4_7}. 

\section{Signal and Orthogonal Bases}\label{sec4_2}

For reasons that will become clear in Section \ref{sec4_3}, it is convenient to represent a signal in a way that allows decoupling the contribution of every time instance. In other words, we seek to define a signal with respect to a very localized \emph{basis} that allows for \emph{sampling} the signal at a given time. With such a (unitary) basis, the sampling process can be done by direct comparison: for example, we say that a continuous signal $u(t)$ has $u(3)=2$ because at time $t=3$ this signal equals two times the element of the basis that is unitary (in a sense to be defined) at $t=3$ and zero elsewhere. Mathematically, this `comparison' process is a \emph{correlation}, the signal processing equivalent of the \emph{inner product}. This notion is more easily introduced for discrete signals, considered in \ref{sec4_2_2}. Continuous signal are treated in \ref{sec4_1_1}.

\subsection{Discrete Signals}\label{sec4_1_1}

Consider a discrete signal of finite duration $u[k]=[1,3,0,4,5]^T$, i.e.\footnote{Note that we will here use a `Python-like' indexing, that is starting from $k=0$ rather than $k=1$.} $k\in[0,\dots 4]$ which can be arranged as a column vector $\mathbf{u}\in\mathbb{R}^{5 \times 1}$. In many applications, a finite duration signal is assumed to be a special case of an infinite duration signal which is zero outside the available points. This is a common practice referred to as \emph{zero-padding} that will be further discussed in Chapter 8. This signal can be written as a linear combination of shifted unitary impulses:

\begin{equation}
	\label{eq1}
	\mathbf{u}=1\, \begin{bmatrix}1 \\0\\0\\0\\0 \end{bmatrix}_{\delta[l-0]} +3\, \begin{bmatrix}0 \\1\\0\\0\\0 \end{bmatrix}_{\delta[l-1]}+0\, \begin{bmatrix}0 \\0\\1\\0\\0 \end{bmatrix}_{\delta[l-2]}+4\, \begin{bmatrix}0 \\0\\0\\1\\0 \end{bmatrix}_{\delta[l-3]}+5\, \begin{bmatrix}0 \\0\\0\\0\\1 \end{bmatrix}_{\delta[l-4]}\,.
\end{equation} 

The elementary basis to describe such a signal is thus the set:

\begin{equation}
	\label{eq2}
	\bm{\delta}_k= \delta[l-k] = \begin{cases} 0 &\mbox{if }  l\neq k \\
		1 & \mbox{if } l= k \end{cases} \quad  \mbox{and}\quad \sum^{l=\infty}_{l=-\infty} \delta[l-k] =1\,.
\end{equation}

Note that two notations are introduced. $\bm{\delta}_k$ is a \emph{vector} of the same size of $\mathbf{u}$ which is zero expect at $l=k$, where it is equal to one; $\delta[l-k]$ is a \emph{sequence} of numbers collecting the same information. In this sequence, $k$ is the index spanning the position of the impulse while the index $l$ spans the time (shift) domain. Infinite duration signals are vectors of infinite length. For two such signals $a[k], b[k] $ or vectors $\mathbf{a},\mathbf{b}$, the inner product for the `comparison procedure' is

\begin{equation}
	\label{eq3}
	\langle \mathbf{a},\mathbf{b}\rangle= \bm {b}^{\dag}\bm{a}=\sum^{n_t-1}_{k=0} a[k] \,\overline{b}[k]\,.
\end{equation} where $b^{\dag}=\overline{b}^T$is the Hermitian transpose, with the superscript $^T$ denoting transposition and the over-line denoting conjugation. With such a basis, the value of the signal at a specific location $u[k]=\mathbf{u}_k$ can be written as

\begin{equation}
	\label{eq4}
	\boxed{
		\mathbf{u}_k=u[k]=\langle \mathbf{u},\bm{\delta}_k\rangle_l \quad \mbox{or}\quad u[k]=\sum^{\infty}_{l=-\infty} u[l]\delta[k-l]
	}\,.
\end{equation}

The operation on the left is a \emph{correlation}: for a given location of the impulse $k$, the inner product is performed over the index $l$ spanning the entire length of the signal and the result is a scalar -- the signal's value at the index $k$. The operation on the right is a \emph{discrete convolution} and the result is a signal: the entire set of shifts will be spanned. We shall return to the algebra of this operation in Chapter 8. 

Note the flipping of the indices $\delta[l-k]$ in \eqref{eq2} to $\delta[k-l]$ in \eqref{eq4}. In the first case, the location of the impulse $k$ is \emph{fixed} and $l$ spans the vector entries; in the second case, within the summation, the time domain $k$ is fixed and $l$ spans the possible locations of the delta functions\footnote{This distinction is irrelevant for a symmetric function such as $\bm{\delta}$, but it is essential in the general case: if $\delta[k-l]$ is replace by $\delta[l-k]$ in \eqref{eq4}, the operation is called cross-correlation.}.

The inner product in a vector space is the fundamental operation that allows for a rigorous definition of intuitive geometrical notions such as the length of a vector and the angle between two vectors. The length ($l^2$ norm) of a vector $||\mathbf{a}||$ and the cosine of the angle $\beta$ between two vectors $\mathbf{a}$ and $\mathbf{b}$ of equal size are defined respectively as

\begin{equation}
	\label{eq5}
	||\mathbf{a}||=\sqrt{\langle \mathbf{a},\mathbf{b}\rangle}=\sqrt{\mathbf{b}^{\dag}\mathbf{a}} \quad \mbox{and} \quad \cos(\beta)=\frac{\langle \mathbf{a},\mathbf{b}\rangle}{||\mathbf{a}||||\mathbf{b}||}
\end{equation}

In signal processing, the first quantity is the root of the signal's \emph{energy}\footnote{Note that the notion of energy is used in signal processing for the square of a signal independently of whether this is actually linked to physical energy.}, defined as $\mathcal{E}\{y\}=||\mathbf{y}||^2$ while the second is the \emph{normalized correlation} between two signals.

\begin{wrapfigure}{r}{0.49\textwidth}
	\begin{center}
		\includegraphics[width=0.2\textwidth]{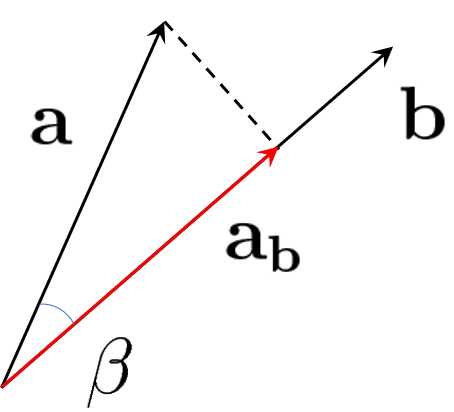}
	\end{center}
	\caption{Projections and inner products.}
	\label{Fig2}
\end{wrapfigure}

Signals with finite energy are \emph{square-summable}. If $\cos \beta=0$, two vectors are \emph{orthogonal} and two signals are \emph{uncorrelated}; if $\cos \beta=1$ two vectors are \emph{aligned} and two signals are \emph{perfectly correlated}.

Notice that the projection of a vector $\mathbf{a}$ onto a vector $\mathbf{b}$ (see Figure \ref{Fig2}) is 

\begin{equation}
	\label{eq6}
	\mathbf{a}_b=||\mathbf{a}|| \cos(\beta)=\frac{\langle \mathbf{a},\mathbf{b}\rangle}{||\mathbf{b}||}\,.
\end{equation}

Hence if $\mathbf{b}$ is a basis vector of unitary length, the inner product equals the projection. When this occurs, as in most of the cases presented in what follows, the notion of \emph{inner product} and \emph{projection} are used interchangeably.

The basis of shifted impulses has a very special property: it is \emph{orthonormal}. This means that the inner product of two basis elements (in this case the shifted delta functions) gives zero unless the same basis element is considered, in which case we recover its unitary norm (energy). We return to this property in Chapter \ref{chapter8}. Before moving to continuous signals, it is worth introducing another important signal that is linked to delta functions, namely the Heaviside step function. This is defined as 

\begin{equation}
	\label{eq7}
	u_{S}[k-l]=\begin{cases} 1 &\mbox{if }  k\geq l \\
		0 & \mbox{if } l< k \end{cases} \, \quad \mbox{i.e.} \longrightarrow u_{S}[k-l]=\sum^{k}_{r=-\infty} \delta [r-l]
\end{equation}

The difference between two shifted step functions generates a box function $u_{\mathcal{B}a,b}$ which is unitary in the range $k=[a,b-1]$ and is zero outside. Tthe difference between two step functions shifted by a single step is a delta function, i.e. $\delta[k-l]=u_{S}[k-l]-u_{S}[k-l-1]$. In discrete setting, this is equivalent to a differentiation. Hence, \emph{the delta function is the derivative of the step function} and the summation in \eqref{eq7} shows that \emph{the step function is an integral of the delta functions}. The reader should close this subsection with an exercise on some distinctive features of discrete signals. 
\medskip

\begin{tcolorbox}[breakable, opacityframe=.1, title=Exercise 1: Continuous vs Discrete]
	
	Consider the continuous signal $y(t)=e^{-0.2 t} \,\cos(2\pi t)$ where $t$ is in seconds. Assuming that the signal is sampled at $f=100 \mbox{Hz}$, and that $n_t=320$ points are collected from $t=0$, what is the resulting discrete signal? Plot this signal and comment on its link with the continuous one.
	
\end{tcolorbox}

\subsection{Continuous Signals}\label{sec4_2_2}

The extension of the previously introduced notion to continuous signals brings several complications, a detailed resolution of which is out of the scope of this chapter. Interested readers are referred to \cite{kaiser2010friendly} for a gentle introduction. The main difficulty arise from the need to define an inner product space which can generalize \eqref{eq3} for functions while allowing for a basis of impulses, i.e. functions that are zero at all but one point. This generalization is provided by the Hilbert space, within which the inner product of two complex-valued functions $a(t), b(t)\in \mathbb{C}$ is

\begin{equation}
	\label{eq8}
	\langle a,b\rangle=\int^{\infty}_{-\infty} a(t)\,\overline{b}(t)\,dt.
\end{equation}

The notion of energy (norm), correlation and projection in \eqref{eq5}-\eqref{eq6} extends to continuous signals using the inner product\footnote{Note that the upper case is used for the norms of a continuous function, i.e. $||a(t)||_2=\sqrt{\langle a,a \rangle}$, with the inner product in \eqref{eq8}, is the $L^2$ norm of the function $a(t)$.} in \eqref{eq8}. The generalization of \eqref{eq4} requires the definition of a continuous delta function $\delta(t)$ with the same properties as $\delta[k]$: it is nonzero only at a given time and  it is absolutely integrable. This can be constructed as the limit of a suitably chosen function having unity area over an infinitesimal time interval. An example is a normalized Gaussian $G(t,\mu,\sigma)$ which has unitary integral regardless of its standard deviation $\sigma$: 

\begin{equation}
	\label{eq9}
	\int^{\infty}_{-\infty}G(t,\mu=0,\sigma)dt=\frac{1}{\sqrt{2\pi\sigma^2}}\int^{\infty}_{-\infty} e^{\frac{-t^2}{2 \sigma^2}} \,dt=1\,.
\end{equation}

This functions gets narrower and taller, as $\sigma\rightarrow 0$, to the point in which it becomes infinite at $t=0$ and null everywhere else, while still having unitary area. This is the definition of continuous Dirac delta function, which in its shifted form is:

\begin{equation}
	\label{eq10}
	\delta(t-\tau)=\begin{cases} 0 &\mbox{if }  t\neq \tau \\
		\infty & \mbox{if } t=\tau \end{cases} \quad \mbox{and} \quad  \int^{+\infty}_{-\infty} \delta(t-\tau) dt=1
\end{equation}

This function is not an ordinary one, as its integration poses several technical difficulties. Without entering into details of measure and distribution theory \cite{Richards1990}, we shall accept this as a \emph{generalized function} that serves well our purpose of sampling a continuous signals. From the definition in \eqref{eq10}, it is easy to derive the sifting (sampling property):

\begin{equation}
	\label{eq11}
	\begin{split}
		\boxed{
			y(t)=\langle y(\tau), \delta (t-\tau) \rangle=\int^{\infty}_{-\infty} y(\tau)\delta(t-\tau)d\tau}=\int^{\infty}_{-\infty} y(t)\delta(t-\tau)d\tau\\=
		y(t)\,\int^{\infty}_{-\infty} \delta(t-\tau)d\tau=y(t)
	\end{split}
\end{equation}

The equivalence $y(\tau)=y(t)$ in the integral results from the product $y(\tau)\delta(t-\tau)$ being null everywhere but at $t$; then, in the last step it is possible to move $y(t)$ outside the integral as this is independent from the integration domain $\tau$.

As for the discrete case, it is interesting to introduce the unitary step function and its link with the delta function as:

\begin{equation}
	\label{eq12}
	u_{S}(t-t_0)=\begin{cases} 1 &\mbox{if }  t> t_0 \\
		0 & \mbox{if } t< t_0 \end{cases} \, \quad \mbox{i.e.} \longrightarrow u_{S}(t-t_0)=\int^{t}_{-\infty} \delta (\tau-t_0) d\tau
\end{equation} Notice that the step function is not defined at $t=t_0$. This definition is the continuous analogue of \eqref{eq7}. To show that the delta function is the derivative of the step function we must introduce the notion of \emph{generalized derivative}. For a continuous signal $u(t)$, denoting as $u'$ and $u^{(n)}$ its first and $n^{th}$ derivative, integration by part using an appropriate test function $\xi(t)$ gives

\begin{equation}
	\label{eq13}\int^{\infty}_{-\infty} \xi(t) u^{(n)}(t) dt = (-1)^{n}\int^{\infty}_{-\infty}\xi^{(n)} u(t) dt
\end{equation} where the test function $\xi(t)$ is assumed to be continuous and differentiable at least $n$ times and is such that $\xi(t)\rightarrow 0$ for $t\rightarrow \pm \infty$. The first derivative of $u_{\mathcal{S}}$ is 

\begin{equation}
	\label{eq14}
	\begin{split}
		\int^{\infty}_{-\infty} \xi(t) \boxed{u'_{\mathcal{S}}(t-t_0)} dt&=-\int^{\infty}_{-\infty} \xi'(t) u_{\mathcal{S}}(t-t_0) dt=-\int^{\infty}_{t_0} \xi^{'}(t)=\\&=\xi(t_0)-\xi(\infty)=\xi(t_0)=\int^{\infty}_{-\infty} \xi(t) \boxed{\delta(t-t_0)} dt
	\end{split}
\end{equation}

The first equality results from direct application of \eqref{eq13}; the second from $u_{\mathcal{S}}(t-t_0)$ being zero in $t<t_0$. After integration, the sifting property of the delta function \eqref{eq11} is used, and the result is obtained by equivalence of last step with the first, holding on the fact that $\xi(t)$ is arbitrary.

\section{Convolutions and Eigenfunctions}\label{sec4_3}

In a \emph{linear system}, the input-output relation satisfy \emph{Homogeneity} and \emph{Superposition}.
For continuous ($\mathcal{S}_c\{u(t)\}=y(t)$) and discrete ($\mathcal{S}_d\{u[k]\}=y[k]$) systems, these set
\bigskip

\textbf{Homegeneity}:
\begin{equation}
	\label{eq15}
	\begin{split}
		\mathcal{S}_c\{a\,u(t)\}=a\mathcal{S}_c\{\,u(t)\}=a y(t)\quad \forall a \in \mathbb{C}\\
		\mathcal{S}_d\{a\,u[k]\}=a\mathcal{S}_d\{\,u[k]\}=a y[k]\quad \forall a \in \mathbb{C}
	\end{split}
\end{equation}

\textbf{Superposition}:
\begin{equation}
	\label{eq16}
	\begin{split}
		&\,\,\mathcal{S}_c\Biggl\{\int^{\infty}_{-\infty}u(t)\Biggr \}dt =\int^{\infty}_{-\infty}\mathcal{S}_d\{u(t) \} dt=\int^{\infty}_{-\infty} y(t) dt\\
		&\,\,\mathcal{S}_d\Biggl\{\sum^{N}_{n=1}u_n[k]\Biggr \}=\sum^{N}_{n=1}\mathcal{S}_d\{u_n[k] \}=\sum^{N}_{n=1} y_n[k]
	\end{split}
\end{equation} where we considered a finite summation of $N$ inputs for the discrete case and an infinite summation of infinitesimally close inputs for the continuous one. Combining these properties, we see that a linear combination of inputs results in the same linear combination of outputs:

\begin{equation}
	\label{eq17}
	\begin{split}
		&\,\,\mathcal{S}_c\Biggl\{\int^{\infty}_{-\infty} a(\tau) u(t,\tau)\Biggr \}d\tau =\int^{\infty}_{-\infty}a(\tau)\mathcal{S}_d\{ u(t,\tau) \} d\tau=\int^{\infty}_{-\infty} a(\tau)y(t,\tau) d\tau\\
		&\,\,\mathcal{S}_d\Biggl\{\sum^{N}_{n=1}a_n u_n[k]\Biggr \}=\sum^{N}_{n=1}a_n \mathcal{S}_d\{u_n[k] \}=\sum^{N}_{n=1} a_n y_n[k]\\
	\end{split}
\end{equation}

A system is \emph{time-invariant} (or \emph{translation-invariant}, if time is replaced by space) if the response to the input does not change over time (or space), i.e.:

\begin{equation}
	\label{eq18}
	\begin{split}
		& \,\, \mbox{if} \, \mathcal{S}_c \{ u(t)\}=y(t) \quad \mbox{then} \quad\mathcal{S}_c \{ u(t-t_0)\}=y(t-t_0)\,\,\,\forall t_0\in \mathbb{R}\\
		& \,\, \mbox{if} \, \mathcal{S}_d \{ u[k]\}=y[k] \quad \mbox{then}\quad\mathcal{S}_d \{ u[k-k_0]\}=y[k-k_0] \,\,\,\forall k_0\in \mathbb{Z}
	\end{split}
\end{equation}

A system is Linear Time/Translation Invariant (LTI) if it is both linear \emph{and} time/translation invariant. These two properties, combined with the signal representations in \eqref{eq10}-\eqref{eq11} make the analysis of LTI systems particularly simple: the response to \emph{any} input can be fully characterized from the response to a \emph{single impulse}.

Defining as $h(t)=\mathcal{S}_c\{\delta(t) \}$ the impulse response of a continuous system, the response to any input is

\begin{equation}
	\label{eq19}
	\begin{split}
		&{y(t)}=\mathcal{S}_c\{ u(t)\}=\mathcal{S}_c\Biggl \{ \int^{\infty}_{-\infty} u(t)\delta(t-\tau) d\tau\Biggr\}=\int^{\infty}_{-\infty} u(t)\mathcal{S}_c \{\delta(t-\tau)\} d\tau\rightarrow\\
		\rightarrow&\boxed{y(t)=\int^{\infty}_{-\infty} u(\tau) \,h(t-\tau) d\tau =\int^{\infty}_{-\infty} h(\tau) \,u(t-\tau) d\tau}
	\end{split}
\end{equation}

The last two integrals are equivalent forms of the convolution integral, hinging on its commutative property.
Similarly, defining $h[k]=\mathcal{S}_d\{\delta[k] \}$ the impulse responses of a discrete system, the response to any input is

\begin{equation}
	\label{eq20}
	\begin{split}
		y[k]&=\mathcal{S}_d\{ u[k]\}=\mathcal{S}_d\Biggl \{ \sum^{\infty}_{l=-\infty} u[k]\delta[k-l]\Biggr\}=\sum^{\infty}_{l=-\infty} u[k]\mathcal{S}_d \{\delta[k-l]\}\rightarrow\\
		\rightarrow&\boxed{y[k]=\sum^{\infty}_{l=-\infty} u[l] \,h[k-l] =\sum^{\infty}_{l=-\infty} h[k] \,u[k-l]}\,,
	\end{split}
\end{equation} having introduced the discrete convolution and its commutative property. 

Note that the system response obtained via \eqref{eq19} or \eqref{eq20} is independent from the initial state of a system. Such response is thus a \emph{particular} solution, i.e. after the transitory from the initial condition vanishes.

A special case is produced when the input is an exponential of the form $u(t)=e^{s t}$ for continuous systems and $u[k]=z^{k}$ for the discrete ones, with $s,z \in \mathbb{C}$. In what follows, we write $s$ in a Cartesian form as $s=\sigma +\mathrm{j} \omega$ and $z$ in a polar form as $z=\rho e^{\mathrm{j}\,\theta}$; the convenience in this is evident in Section \ref{sec4_6}. For the moment, note that these two variables are the continuous and discrete \emph{complex frequencies}.

In continuous systems, from \eqref{eq19}, the output is

\begin{equation}
	\label{eq21}
	y(t)=\mathcal{S}_c\{e^{s t}\}=\int^{\infty}_{-\infty} h(\tau) \,e^{(t-\tau) s} d\tau = \Biggl(\int^{\infty}_{-\infty} h(\tau) \,e^{-\tau s} d\tau \Biggl)\,e^{ s t} =\lambda_c(s)e^{ s t}
\end{equation}

In discrete systems, from in \eqref{eq20}, the output is

\begin{equation}
	\label{eq22}
	y[k]=\mathcal{S}_d\{z^{k}\}=\sum^{\infty}_{l=-\infty} h[l] \,z^{k-l} = \Biggl(\sum^{\infty}_{l=-\infty} h[l] \,z^{-l} \Biggl)\,z^{k} =\lambda_d(z) z^{k}
\end{equation}

In both cases, this result shows that LTI system responds to a complex exponential with \emph{the same input} multiplied by a complex number ($\lambda_d$ or $\lambda_c$). This number solely depends on the complex frequencies $z \in \mathbb{C}$ and $s\in \mathbb{C}$.
Therefore, these special inputs functions are \emph{eigenfunctions} of the LTI operators and their complex eigenvalues are: 

\begin{subequations}
	\label{eq23}
	\begin{equation}
		\lambda_c(s)=H(s)=\int^{\infty}_{-\infty} h(\tau) \,e^{-\tau s} d\tau
	\end{equation}
	\begin{equation}
		\lambda_d(z)=H(z)=\sum^{\infty}_{l=-\infty} h[l] \,z^{-l}
	\end{equation}
\end{subequations}

These are respectively the Laplace transform and the Z-Transform of the impulse response. These are the \emph{transfer functions} of the LTI systems and link input and output in the \emph{complex} frequency domain. For time varying or nonlinear systems, the notion of transfer function is not useful. Finally, note that in discrete systems, the transfer function is a continuous function of $z$.

\bigskip

\begin{tcolorbox}[breakable, opacityframe=.1, title=Exercise 2: Convolution with Impulse Responses]
	
	To identify a SISO system, a step function test is carried out. A regression analysis reveals that the response of the system to a step $u_{\mathcal{S}}$ is well described by the function 
	\begin{equation}
		\label{eq24}
		y_{\mathcal{S}}=\frac{1}{5}-\frac{1}{5}e^{-t}\bigl(\cos(2\,t)+\frac{1}{2}\sin(2\,t)\bigr)
	\end{equation}
	Compute the response of the system to an exponential input $u(t)=e^{-t}$ in the continuous and the discrete domains, assuming that the system can be sampled at $f_s=10 \mbox{Hz}$ or at $f_s=3 \mbox{Hz}$ for $T=5 s$.

\end{tcolorbox}

\section{Causal and Stable Systems}\label{sec4_4}

In most applications of interest, signals are assumed to be null at time $t<0$. This is important for the impulse response of systems that are \emph{causal}, in which the impulse response is $h(t)=0$ for $t<0$ and $h[k]=0$ for $k<0$. This means that no output can be produced \emph{before} the input, and hence the system is not anticipatory. Models of physical systems and \emph{online} data processing must be causal. On the other hand, many data processing schemes operating \emph{offline} are \emph{not causal} (e.g., zero-phase filters described in section \ref{sec4_7}). 

As anticipated in the previous exercise, the convolution integral and summations in \eqref{eq19}- \eqref{eq20} for \emph{casual} systems become

\begin{subequations}
	\label{eq25}
	\begin{equation}
		y(t)=\mathcal{S}_c\{ u(t)\}=\int^{t}_{0} h(\tau) \,u(t-\tau) d\tau=\int^{t}_{0} u(\tau) \,h(t-\tau) d\tau
	\end{equation}
	\begin{equation}
		y[k]=\mathcal{S}_d\{ u[k]\}=\sum^{k}_{l=0} h[l] \,u[k-l]=\sum^{k}_{l=0} u[l] \,h[k-l]\,.
	\end{equation}
\end{subequations} 

The upper limit is replaced by $t$ or $k$ since signals and impulse responses are null for $\tau>t$ or $k>l$ and the lower one is replaced by $0$ since both are null for $t<0$ or $l<0$. The continuous and discrete transfer functions in \eqref{eq23} become 

\begin{equation}
	\label{eq26}
	\lambda_c(s)=H(s)=\int^{\infty}_{0^{-}} h(\tau) \,e^{-\tau s} d\tau \quad \mbox{and} \quad \lambda_d(z)=H(z)=\sum^{\infty}_{n=0} h[n] \,z^{-n}
\end{equation}

In the continuous case, $t=0$ \emph{must} be included in the integration; hence the lower bound is tuned to accommodate for any peculiarity occurring at $t=0$ (notably an impulse). Nevertheless, to avoid the extra notational burden, the minus subscript in the lower bound is dropped in what follows.

Finally, another important class of interest is that of \emph{stable} systems. The stability analysis of complex systems is a broad topic (see Chapters 10 and 13). Here, we limit the focus to bounded-input/bounded output (BIBO) stability. A system is BIBO stable if its response to \emph{any} bounded input is a bounded output. This requires that the impulse response of continuous and discrete signals satisfy:

\begin{equation}
	\label{eq27}
	\int^{\infty}_{0} |h(t)|dt<\infty\, \quad \mbox{and} \quad \sum^{\infty}_{k=0} |h[k]|<\infty\,.
\end{equation}

A different notion is that of \emph{asymptotic} stability, which is related to the \emph{internal} stability of a system.
A system is \emph{asymptotically stable} if every initial state, in the absence of inputs, produces a bounded response that converges to zero. In an LTI system, asymptotic stability implies BIBO stability, but the reverse is not true.

\section{LTI Systems in their Eigenspace}\label{sec4_5}

In the previous section, we have seen that complex exponentials are eigenfunctions of LTI systems.
Great insights on a system behavior can be obtained by projecting their input-output relation onto the system's eigenfunctions. The projection of signals into complex exponentials leads to the Laplace transform in the continuous domain and the Z transform in the discrete domain. This section is divided into four subsections. We start with some definitions.

\subsection{Laplace and Z Transforms}\label{sec4_5_1}

\textbf{Laplace Transforms.} Given a continuous signal and causal signal $u(t)$, the Laplace transforms are

\begin{equation}
	\label{eq28}
	U(s)=\mathcal{L}\{u(t)\}=\int^{\infty}_{-\infty} u(t) e^{-st}dt=\int^{\infty}_{0} u(t) e^{-st}dt\,.
\end{equation}

The first integral is the \emph{bilateral} transform; the second is the \emph{unilateral} transform. As we here focus on causal signals (i.e. $u(t<0)=0$), these are identical. Nevertheless, these are different tools required for different purposes: the first is suitable for infinite duration signals, for which it can be linked to the Fourier Transform (Section \ref{sec4_6}); the second is developed for solving initial value problems, as it naturally handles initial conditions. 

These integral converge, and hence the Laplace transforms exist, if $u(t)e^{-st}\rightarrow 0$ for $t\rightarrow \pm \infty$ (only $t\rightarrow +\infty$ for the \emph{unilateral}). This requires that the signal is of \emph{exponential order}, i.e. grows more slowly than a multiple of some exponential: $|u(t)|\leq M e^{\alpha t}$. If this is the case, the range of values $\mathbb{R}\{s\}>\alpha$ is the \emph{region of convergence} (ROC) of the transform. 

The reader is referred to \cite{Beerends2003}, \cite{Wang2009}, \cite{Hsu2013} for a review of all the properties of the Laplace transform; we here focus on the key operations enabled by this powerful tool and we omit formulation of the inverse Laplace transform, as it requires notions of complex variables theory that are out of the scope of this chapter. The Python script \textsc{Ex3.py} to solve Exercise 3 provides the commands to compute both the transform and its inverse using the Python library \textsc{SymPy}\footnote{see \url{https://www.sympy.org/en/index.html}.}.

The key property of interest in this chapter is that of time derivation, which can be easily demonstrated using integration by parts. The bilateral ($\mathcal{L}^{b}$) and unilateral ($\mathcal{L}^{u}$) transforms of a time derivative are

\begin{equation}
	\label{eq29}
	\mathcal{L}^{b}\{u'(t)\}=s U(s) \quad \mbox \quad \mathcal{L}^{u}\{u'(t)\}=s U(s)-u(0)\,.
\end{equation} 

That is \emph{differentiation in the time domain corresponds to multiplication by $s$ in the frequency domain}; similarly, one can show that \emph{integration in the time domain corresponds to division by $s$}. Notice that no distinction between $\mathcal{L}^{b}$ and $\mathcal{L}^{a}$ is needed for a system initially at rest and the \emph{bilateral transform cannot handle initial conditions}. 

Finally, compare the inner product in \eqref{eq28} with \eqref{eq8}, taking $a(t)=u(t)$. Note that the Laplace transform is a projection of the signal $u(t)$ onto an exponential basis $b_{\mathcal{L}}(t,s)=e^{-\sigma +\mathrm{j}\omega}$. It is left as an exercise to show that the Laplace basis is not orthogonal, unless $\sigma=0$. That is the basis of the continuous Fourier transform.  

\bigskip
\textbf{Z Transform.} Given a discrete and causal signal $u[k]$, the Z transforms are: 

\begin{equation}
	\label{eq30}
	U(z)=\mathcal{Z}\{u[k]\}=\sum^{+\infty}_{k=-\infty} u[k] z^{-k}=\sum^{\infty}_{k=0} u[k] z^{-k}\,.
\end{equation}

The same distinction on \emph{bilateral} or \emph{unilateral} transforms, as well as their equivalence for causal signals, apply to the discrete transform. Observe that the Z transform is a \emph{continuous function} of $z$. These summation converge, and hence the transforms exists, if $|u[k] z^{-k}|<1$. The domain within which this occurs is the \emph{region of convergence} (ROC) of the transforms, and this is typically within a domain such that $|z|>\alpha$. As in the previous section, we avoid a review of all the properties of this transform, and the formulation of its inverse (see \cite{Wang2009,Hayes2011}).

The Z-transform equivalent of the time differentiation and integration properties of the Laplace transform are the time-shifting properties:

\begin{equation}
	\label{eq31}
	\mathcal{Z}\{u[k-1]\}=z^{-1}\, U(z) \quad \mbox{and} \quad \mathcal{Z}\{u[k+1]\}=z\, U(z)\,.
\end{equation}

Finally, compare \eqref{eq30} with \eqref{eq3} taking $a[k]=u[k]$ to see that the Z-transform is a projection of the signal onto a basis of powers $b_{\mathcal{Z}}[k]=\rho^{-k}e^{\mathrm{j}\theta k}$. This basis is not orthogonal unless $\rho=1$. That is the basis of the discrete Fourier transform. 

The Z-transform is not suited for compression purposes, as the dimensionality of the problem in the transformed domain is \emph{increased}: the original signal has a (finite) dimension $t_k$, but the same information in the frequency domain is mapped onto a \emph{continuous complex plane} $z$.

\subsection{Discrete and Continuous Frequencies}\label{sec4_5_2}

If a discrete signal is obtained by sampling a continuous one, the link between these two transforms reveals the important impact of the sampling process on the frequency domain. The discretization $t_k=k\Delta t$ creates a point-wise equivalence such that any signal (causal or not) can be equivalently written as

\begin{equation}
	u[k]=\sum^{\infty}_{l=-\infty} u[l]\delta[k-l]\longleftrightarrow u(t_k)=\sum^{\infty}_{l=-\infty} u(t_l)\delta(t_k-l\Delta t)\,.
\end{equation}

With $e^{-s \,t_k} =e^{-s \,k \Delta t}=z^k$, the Laplace transform of this discrete signal is

\begin{equation}
	\label{e33}
	\begin{split}
		\mathcal{L}\{u[k]\}&=\int^{\infty}_{-\infty} \sum^{\infty}_{l=-\infty} u[l] \delta(t-l\Delta t) e^{-s t} dt=\\
		&=\int^{\infty}_{-\infty} \sum^{\infty}_{l=-\infty} u[l] \delta(t-l\Delta t) z^{-l} dt=\\
		&=\sum^{\infty}_{l=-\infty} u[l] z^{-l} \int^{\infty}_{-\infty} \delta(t-l\Delta t)dt=\mathcal{Z}\{u[k]\}
	\end{split}
\end{equation}

The equivalence of these transforms relies on the change of variables $z=e^{s \Delta t}$. This maps $s$ onto the $z$ while preserving angles in the two domains: such mapping is a \emph{conformal mapping}. Introducing $s=\sigma+\mathrm{j}\omega$ and $z=\rho e^{\mathrm{j}\theta}$, shows that $\rho=e^{\sigma \Delta t}$ and $\theta=\omega \Delta t$.

\begin{figure}[htbp]
	\centering
	\includegraphics[keepaspectratio=true,width=0.81 \columnwidth]{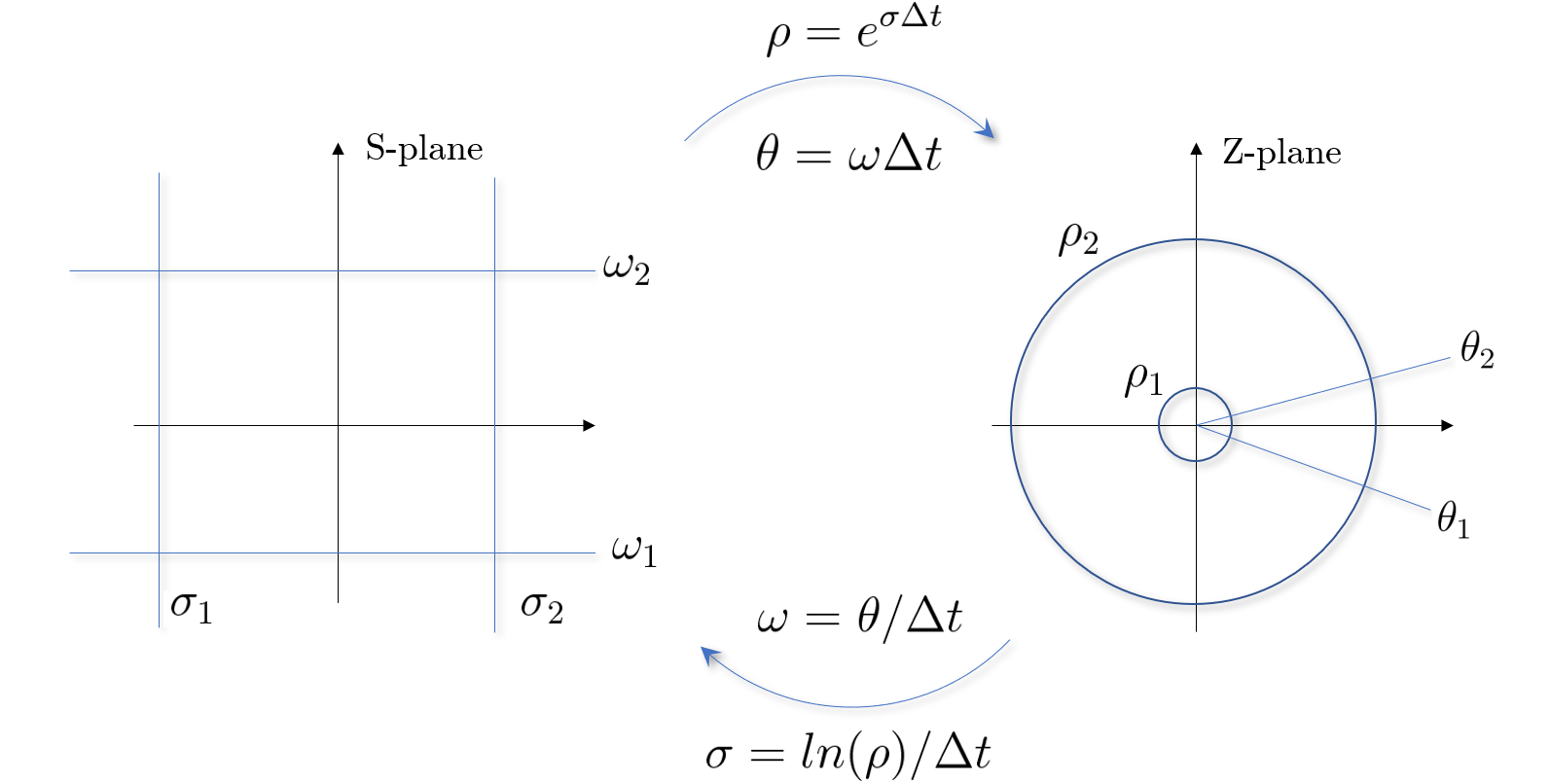}
	\caption[Examples of Discrete Signals]{The conformal mapping linking the Laplace Transform and the Z transform.}
	\label{fig5}
\end{figure}

Figure \ref{fig5} shows several important features of this mapping. Observe that the imaginary axis $\omega$ in the continuous domain is mapped onto an angular coordinate $\theta\in [-\pi ,\pi]$ or $\theta\in [0,2\pi]$. This offers yet another way of introducing the notion of \emph{aliasing} and the Nyquist-Shannon sampling theorem (encountered in the exercise 1), arising from the fact that the frequency domain of a digital signal is periodic. Another important observation is that the left side of the complex domain $s=\sigma<0$ is mapped \emph{inside} the unit circle $\rho=1$, while the axis $s=\sigma=0$ is mapped \emph{on} the unit circle.

\subsection{The Convolution Theorem}\label{sec4_5_3}

While the convolution integral and summations in \eqref{eq19} and \eqref{eq20} provide the input-output relation in the time domain using the impulse response, further insights can be obtained by analyzing this relation in the frequency domain. We here give, without proofs\footnote{See \cite{Wang2009} for more details.}, one of the most important results of signal processing, known as the convolution theorem: \emph{a convolution in the time domain is a multiplication in the frequency domain}. For continuous and discrete systems, this means:

\begin{subequations}
	\label{eq34}
	\begin{equation}
		Y(s):=\mathcal{L}\{y(t)\}=\mathcal{L}\Biggl \{\int^{t}_{0} h(\tau) \,u(t-\tau) d\tau\Biggr \}=H(s) U(s)
	\end{equation}
	\begin{equation}
		Y(k):=\mathcal{Z}\{y[k]\}=\mathcal{Z}\Biggl \{\sum^{k}_{l=0} h[l] \,u[k-l] \Biggr \}=H(z) U(z)\,.
	\end{equation}
\end{subequations}

\subsection{Differential and Difference Equations}\label{sec4_5_4}

Continuous LTI systems can be described in terms of (linear) differential equations while discrete LTI systems can be described in terms of (linear) difference equations. The Laplace and the Z-transform are powerful tools to solve these equations because of the properties in \eqref{eq29} and \eqref{eq31}. In both cases, the equations have constant coefficients and are often acronymized as LCCDE. 

\bigskip
\textbf{Differential Equations.} The general form of the LCCDE of a continuous SISO LTI systems with input $u(t)$ and output $y(t)$ reads

\begin{equation}
	\label{eq35}
	\sum^{N_b}_{n=0}\,a_n\,y^{(n)}(t)=\sum^{N_f}_{n=0}\,b_n\,u^{(n)}(t)\,
\end{equation} where $N_f\geq N_b$ is the order of the system\footnote{The condition $N_f\geq N_b$ is necessary to ensure that the system is \emph{realizable}, that is both stable and causal. More about this in the next footnote.}. The coefficients $a_n$ are called \emph{feedback coefficients}; the coefficients $b_n$ are \emph{feedforward coefficients}.

The LCCDE provides an \emph{implicit} representation of a system since the input-output relation can be revealed only by solving the equation. Introducing the Laplace transform in a LCCDE is an operation similar to the Galerkin projection underpinning Reduced Order Modeling (ROM, see Chapters 1 and 14). Recalling that the Laplace transform is a projection onto the basis $b_{\mathcal{L}}(t)$, \eqref{eq35} leads to:

\begin{equation}
	\label{eq36}
	\begin{split}
		\Bigl \langle \sum^{N_b}_{n=0}a_n y^{(n)}(t), b_{\mathcal{L}}(t)\Bigr \rangle=\Bigl \langle \sum^{N_f}_{n=0}b_n u^{(n)}(t), b_{\mathcal{L}}(t)\Bigr \rangle\rightarrow\\
		\rightarrow \sum^{N_b}_{n=0}a_n \mathcal{L}\Bigl \{y^{(n)}(t)\Bigr \}=\sum^{N_f}_{n=0}b_n \mathcal{L}\Bigl \{u^{(n)}(t)\Bigr \}\rightarrow
		\\\rightarrow Y(s)\Bigl( \sum^{N_b}_{n=0} a_n s^{(n)}\Bigr)=U(s)\Bigl( \sum^{N_f}_{n=0} b_n s^{(n)}\Bigr)\rightarrow \\H(s)=\frac{Y(s)}{U(s)}=\frac{\sum^{N_f}_{n=0} b_n s^{(n)}}{\sum^{N_b}_{n=0} a_n s^{(n)} }=\frac{b_{N_f}}{a_{N_b}}\,\frac{\prod^{N_f}_{n=0} (s-z_n)}{ \prod^{N_b}_{n=0} (s-p_n)}\,,.
	\end{split}
\end{equation} The transfer function of LTI systems is a polynomial rational function of $s$, with the coefficients of the polynomials being the coefficients of the LCCDE. In the factorized form, $z_n$ and $p_n$ are respectively the \emph{zeros} and the \emph{poles} of the system\footnote{A transfer function that has more zeros than poles (i.e. $N_f>N_b$) is said to be \emph{improper}. In this case, $\lim_{s\to\infty} |H(s)|=+\infty$, which violates stability: this implies that at large frequencies, a finite input can produce an infinite output. Moreover, after the polynomial division, the transfer function brings polynomial terms in $s$. The inverse Laplace transform of these are (generalized) derivatives of the delta functions; hence the corresponding impulse response $h(t)=\mathcal{L}^{-1}\{H(s)\}$ violates causality.}. Note that since the coefficients $a_n$, $b_n$ are real, these can either be purely real or appear in complex conjugate pairs. These coefficients have a straightforward connection with the LCCDE, which can immediately be recovered from the transfer function. The zeros $z_n$ are associated to inputs $e^{z_n t}$ in which the transfer function is \emph{null} and thus leads to no output; the poles $p_n$ corresponds to resonances, inputs $e^{p_n t}$ in which the transfer function is infinite and leads to the blow-up of the system. 

The poles are eigenvalues of the matrix  $\mathbf{A}$, advancing a linear system in its state-space representation (see Chapters 10 and 12). 
In a stable system, poles are located in regions of the $s$-plane that are `not-accessible' by any input, that is outside the ROC of the transfer function. Defining the ROC of $H(s)$ as $\mathbb{R}\{s\}>\alpha$, and observing that the poles are by definition outside the ROC, stability is guaranteed if $\alpha=0$, i.e. if \emph{the ROC includes the imaginary axis}. This is equivalent to impose that all the poles are located in the left side of the s-plane, i.e. $\mathbb{R}\{p_n\}<0 $ $\forall n\in[0,N_f]$. This can also be derived from the BIBO stability condition in \eqref{eq27}.

\bigskip
\textbf{Difference Equations.} In the discrete case, the general form of LCCDE associated to SISO LTI systems with input $u[t]$ and output $y[t]$ reads

\begin{equation}
	\label{eq37}
	\sum^{N_b}_{n=0} a_n y[k-n]=\sum^{N_f}_{n=0} b_n u[k-n] , \,\,\mbox{i.e.} \,\, y[k]=\sum^{N_b}_{n=0} b^*_n u[k-n]-\sum^{N_f}_{n=1} a^*_n y[k-n].
\end{equation} The order of the system\footnote{Note that in \eqref{eq37} the restriction $N_b\geq N_f$ is not needed to enforce causality: by construction, the output $y[k]$ only depends on past information.} is $max(N_b,N_f)$. The form on the right plays a fundamental role in filter implementation, time series analysis and system identification and is known as \emph{recursive form} of the differnce equation. Note that the \emph{feedback} and \emph{feedforward} coefficients in the recursive form are simply $a^{*}_n=a_n/a_0$ and $b^{*}_n=b_n/a_0$ respectively, hence the coefficient $b^*_0$ is the \emph{static gain} of the system.

As for the continuous case, projecting \eqref{eq37} onto the $Z$ basis $b_{\mathcal{Z}}(t)$ via Z transform and using \eqref{eq31}, yield the transfer function of a discrete system:

\begin{equation}
	\label{eq38}
	\begin{split}
		\Bigl \langle \sum^{N_b}_{n=0}a_n y[k-n], b_{\mathcal{Z}}[k]\Bigr \rangle=\Bigl \langle \sum^{N_f}_{n=0}b_n u[k-n], b_{\mathcal{Z}}[k]\Bigr \rangle\rightarrow\\
		\rightarrow \sum^{N_b}_{n=0}a_n \mathcal{Z}\Bigl\{ y[k-n]\Bigr \}=\sum^{N_f}_{n=0}b_n \mathcal{Z}\Bigl \{u[k-n]\Bigr \}\rightarrow
		\\\rightarrow Y(z)\Bigl( \sum^{N_b}_{n=0} a_n z^{-n}\Bigr)=U(z)\Bigl( \sum^{N_f}_{n=0} b_n z^{-n}\Bigr)\rightarrow \\H(z)=\frac{Y(z)}{U(z)}=\frac{ \sum^{N_f-1}_{n=0} b_n z^{-n}}{\sum^{N_b-1}_{n=0} a_n z^{-n}}=\frac{b_0}{a_0}z^{N_b-N_f}\,\frac{\prod^{N_f-1}_{n=0} (1-\zeta_n z^{-1})}{ \prod^{N_b-1}_{n=0} (1-\pi_n z^{-1})}
	\end{split}
\end{equation} where $\zeta_n$ and $\pi_n$ are respectively the zeros and poles of the discrete transfer function. Observe that the factored form of the discrete transfer function is usually given in terms of polynomials of $z^{-1}$ rather than $z$. 

The link between zero and poles in continuous and discrete domains is given by the conformal mapping in Figure \ref{fig5}. 
In the absence of inputs, the poles control the evolution of a linear system from its initial condition (i.e. the homogeneous solution of the LCCDE). The Dynamic Mode Decomposition (DMD) introduced in Chapter 7 is a powerful tool to identify the poles $\pi_n$ of a system from data, and to build linear Reduced-Order Models by projecting the data onto the basis of eigenfunctions $z^{\pi_n\,k}$. 

Finally, in analogy with the continuous case, a discrete system is stable if its poles are outside the ROC of the transfer function. Defining the ROC of $H(z)$ as $|z|\geq\alpha$, one sees that this occurs if $\alpha=1$: the ROC include the unit circle and hence all the poles have $|\pi_n|<1$. This can be derived from the BIBO stability condition in \eqref{eq30}.

\bigskip

\begin{tcolorbox}[breakable, opacityframe=.1, title=Exercise 3: Transfer Function Analysis]
	
	Consider the system in Exercise 2. Compute the transfer function and the system output from the frequency domain, then identify the LCCDE governing the input-output relation. Then, assuming that a discrete system is obtained sampling the continuous domain, derive a recursive formula that mimics the input-output link of the continuous system. Test your result for a sampling frequency of $f_s=3Hz$ and $f_s=10Hz$.

\end{tcolorbox}

\section{Application I: Harmonic Analysis and Filters}\label{sec4_6}

BIBO stability guarantees that the output produced by a stationary input is also stationary. It is thus interesting to consider only the portion of the complex planes $s$ and $z$ associated to infinite duration signals, i.e. $s=\mathrm{j}\omega$ and $z=e^{\mathrm{j}\theta}$. These correspond to harmonic eigenvalues of the LTI system, hence lead to an \emph{harmonic response}. From Laplace and Z transforms, we move to continuous and discrete Fourier transforms in section \ref{sec4_6_1}. A system that manipulates the harmonic content of a signal is a filter; these are introduced Section \ref{sec4_6_2} along with their fundamental role in multi-resolution decompositions.

\subsection{From Laplace to Fourier}\label{sec4_6_1}

Consider the bilateral Laplace and Z transform of a signal along the imaginary axis $s=\mathrm{j}\omega$ and the unitary circle $z=e^{\mathrm{j}\theta}$ respectively:

\begin{subequations}
	\begin{equation}
		\label{eq39a}
		U(\mathrm{j}\omega)=\mathcal{F}_C\{u(t)\}=\int^{\infty}_{-\infty} u(t) e^{-\mathrm{j}\omega t}dt
	\end{equation}
	\begin{equation}
		\label{eq39b}	
		U({\mathrm{j}\theta})=\mathcal{F}_D\{u[k]\}=\sum^{+\infty}_{k=-\infty} u[k] e^{-\mathrm{j}\theta\,k}
	\end{equation}
\end{subequations} These are the continuous (CT) and the Discrete (Time) Fourier Transforms (DTFT). Both are continuous functions, with the second being periodic of period $2\pi$ because of the conformal mapping introduced in Figure \ref{fig5}. Comparing these to \eqref{eq28} and \eqref{eq30} shows that the bilateral Laplace and Z transforms are the Fourier transforms of $u(t) e^{-\sigma t}$ and $u[k] \rho^{-k}$. Without these exponentially decaying modulations, the conditions for convergence are more stringent: signals must be absolutely integrable and absolutely summable\footnote{This condition is sufficient but not necessary: some non-square integral functions do admit a Fourier transform. Important examples are the constant function $u(t)=1$ or the step function $u_{\mathcal{S}}(t)$. Moreover, note that the Fourier Transform can be obtained from the Laplace and Z transform, only for signals that are absolutely integrable or summable. For instance, the Laplace transform of $e^{\alpha t}$ with $\alpha>0$, has ROC $s>\alpha$ while the Fourier transform does not exist.}.

The main consequence is that \emph{infinite duration stationary signals do not generally admit a Fourier transform}. This explains why the manipulations of these signals by an LTI system are better investigated in terms of some of their statistical properties, such as \emph{autocorrelation} or \emph{autocovariance}, as illustrated in section \ref{sec4_7}. A special exception are periodic signals, for which \eqref{eq39a} and \eqref{eq39b} lead to Fourier \emph{series}, and the problem of convergence becomes less stringent.

In stable continuous and discrete LTI systems, satisfying \eqref{eq30}, the impulse response always admit Fourier transform: these can be obtained by replacing $s=\mathrm{j}\omega$ and $z=e^{\mathrm{j}\theta}$ in the transfer function. This leads to the \emph{frequency transfer function}, which are complex functions of real numbers\footnote{These are often called \emph{real} frequencies as opposed to the \emph{complex frequencies} $s$ and $z$.} ($\omega$ or $\theta$), customarily represented by plotting $log(|H(x)|)$ and $\mbox{arg} (H(x))$ versus $log(x)$ in a \emph{Bode plot}, with $x=\omega$ or $x=\theta$. The modulus of the frequency transfer function is the \emph{amplitude response}; its argument is the \emph{phase response}. 

If the Fourier transform (or series) exist for both inputs and output, the properties of the Laplace and Z transform applies: the harmonic contents of the output is $Y(\mathrm{j}\omega)=H(\mathrm{j}\omega)U(\mathrm{j}\omega)$ in the continuous domain; $Y(\mathrm{j}\theta)=H(\mathrm{j}\theta)U(\mathrm{j}\theta)$ in the discrete one.

Discrete signals of finite duration $u(t_k)=\mathbf{u}\in\mathbb{R}^{n_t\times 1}=u[k]$, with $k\in[0,n_t-1]$, are usually extended to infinite duration signals assuming periodic boundary conditions. The frequency domain is thus discretized into bins $\theta_n=n\Delta \theta$ with $n\in[0,n_f-1]$ and $\Delta \theta=2\pi/n_F$. The mapping to the continuous frequency domain, from Fig. \ref{fig5}, gives $f_n=2\pi \omega_n=n f_s/n_f$, with $f_s=1/\Delta t$ the sampling frequency. With both time and frequency domain discretized, the Fourier pair are usually written as

\begin{equation}
	\label{eq40}
	U[n]=\frac{1}{\sqrt{n_t}}\sum^{n_t-1}_{k=0} u[k] e^{-2\pi\mathrm{j}\frac{n\,k}{n_f}} \Longleftrightarrow u[k]=\frac{1}{\sqrt{n_t}}\sum^{n_f-1}_{n=0} U[n] e^{2\pi\mathrm{j}\frac{n\,k}{n_f}} \,.
\end{equation}

The equations in \eqref{eq50} are respectively the Discrete Fourier Transform (DFT) and its inverse. Note that the normalization $1/\sqrt{n_t}$ is used for later convenience: we see in Chapter 8 that \eqref{eq50} can be written as matrix multiplications with the columns of the matrix being orthonormal vectors. Finally, if $n_f=n_t$ and $n_t$ is a power of $2$, this multiplication can be performed using the famous FFT (Fast Fourier Transform) algorithm (see \cite{Loan1992}), reducing the computational cost from $n_t^2$ to $n_t log_2(n_t)$. An excellent review of the DFT is provided by \cite{DFT_Smith}.

\bigskip

\begin{tcolorbox}[breakable, opacityframe=.1, title=Exercise 4: Frequency Response Analysis]
	
	Consider the discrete system derived in Exercises 3, but now assume that the static gain is unitary. Compute the frequency transfer function of this system and show that this can be seen as a \emph{low-pass filter}. Study how the frequency response changes if the coefficients $a_1$ or $a_2$ are set to zero. Finally, derive the system that should have the \emph{complementary} transfer function and show its amplitude response. Is this response also complementary? 
\end{tcolorbox}

\subsection{Multiresolution Analysis and Digital Filters}\label{sec4_6_2}
Filters are at the center of most signal processing applications, and the theory behind their design is a vast subject \citep{FIR_Smith}. Among the essential applications discussed in this book are feedback control design (Chapter 10), Multiresolution Analysis (MRA) and wavelet decomposition (Chapter 5) and Multiscale Proper Orthogonal Decomposition (Chapter 8).

In feedback control, a controller manipulates the feedback coefficients of an LCCDE by introducing a control input which is function of the output. Therefore, filter design methods can be used to design the actuation such that the transfer function of the controlled system rejects certain disturbances (see Bode design methods \citep{Distefano2013}).

In MRA, filters are used to decompose signals. While the DFT represents a signal as a linear combination of harmonics, MRA represents it as a combination of frequency bands called \emph{scales}. A packet of similar frequencies can be assembled into bases called \emph{wavelets}, hence the connection to Chapter 5.

The MRA partitions the spectra of signal into $n_M$ scales, each taking a portion of the signal's content in bands $[0,f_1], [f_1,f_2],\cdots [f_{M-1},f_s/2]$. In Chapter 8, these are identified by a \emph{frequency splitting vector} $F_V=[f_1,f_2,f_3\dots f_{M-1}]$.

The MRA of a discrete signal can be written as:

\begin{equation}
	\label{eq41}
	u[k]=\sum^{M}_{m=1} s_{m}[k]=\mathcal{F}^{-1}\Biggl \{\sum^{M}_{m=1} H_m (f_n)U(f_n) \Biggr \}\, \mbox{with} \,\sum^{M}_{m=1} |H_m(f_n)| =1\,,
\end{equation} where $s_m$ is the portion of the signal in the scale $m$, within the frequency range $f_n\in[f_{m-1},f_m]$ and $H_m(f_n)$ is the transfer function of the filter that isolates that portion. Therefore, $|H_m(f_n)|\approx 1$ for $f_n\in[f_{m-1},f_m]$ and $|H_m(f_n)|\approx 0$ otherwise. The assumption on the right enables a lossless decomposition. 

The MRA requires the definition of one low pass filter for the range $[0,f_1]$, one high-pass filter for the range $[f_{M-1},f_s/2]$ and $M-2$ bandpass filters. Because these are complementary, all these filters can be obtained from a set of low pass filters, as described at the end of this section. Therefore, to learn MRA, one should first learn how to construct a low pass filter with a given cut-off frequency $f_c$. 

We now focus on the two main families of filters and the most common design methods. Let us consider a specific example, with $f_s=2kHz$ and $f_c=200 Hz$. In the digital frequency domain, we map the sampling frequency to $\theta_s=2\pi$ and the cut-off to $\theta_c=\pi/5$.  

The transfer function of the ideal low pass filter is

\begin{equation}
	\label{eq42}
	H_{id}(\mathrm{j}\theta)= \begin{cases} e^{-\mathrm{j}\alpha_d \theta} &\mbox{if }  |\theta|\leq \theta_c \\
		0 & \mbox{if } \theta_c<|\theta|\leq\pi \end{cases} 
\end{equation}

This leads to the impulse response 

\begin{equation}
	\label{eq43}
	h_{id}[k]=\frac{1}{2\pi}\int^{\pi}_{-\pi} H(\mathrm{j}\theta) e^{-\mathrm{j}k \theta} d\theta=\frac{\sin[k-\alpha_d]\theta_c}{\pi(k-\alpha_d)}\,.
\end{equation}

\begin{figure}[htbp]
	\centering
	\subfigure
	{
		\includegraphics[keepaspectratio=true,width=0.47 \columnwidth]
		{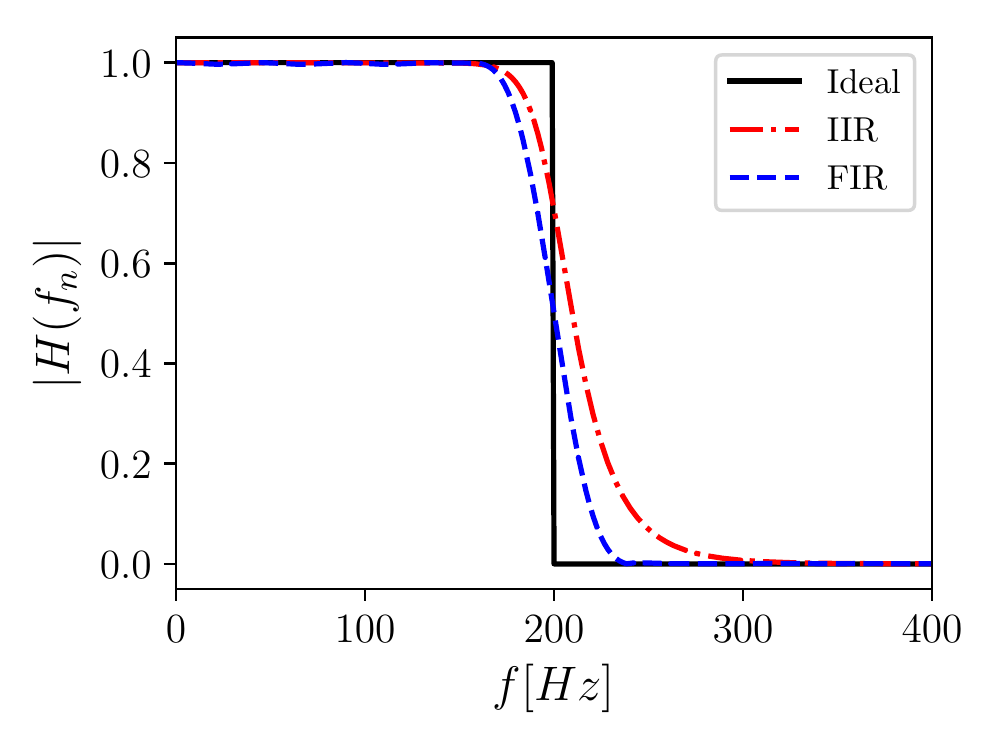}		
	}
	\subfigure
	{
		\includegraphics[keepaspectratio=true,width=0.47 \columnwidth]
		{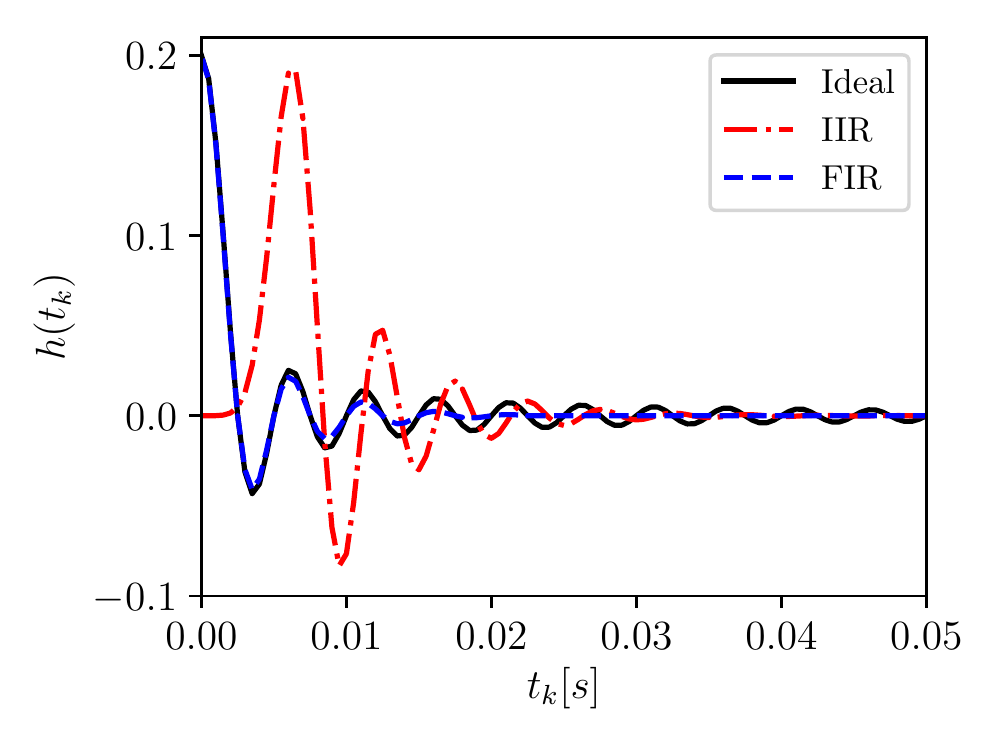}
		
	}
	\vspace{-0.5mm}
	\caption
	{Fig (a): Absolute value of the frequency response of an ideal (black continuous), a IIR (red dash-dotted) and a FIR (blue dashed) low pass filter; Fig (b): impulse responses associated to the frequency responses on the left. In the case of the FIR, the response is shifted backward by $\alpha_d$. }
	\label{fig7}
\end{figure}

The need for a delay $\alpha$ is evident after introducing FIR filters. Notice that having a linear phase $\alpha \theta$ delays the input without distorting its waveform. The modulus of this frequency response and a portion of its transfer function are shown in Figure \ref{fig7}a) with continuous black curves.

Such ideal filter is not realizable in the time domain because its impulse response is not causal ($h[k]\neq0$ for $k<0$) and is not absolutely summable. The ideal constraints must be relaxed. The most popular categories are IIR and FIR filters.

\textbf{1) Infinite Impulse Response Filters (IIR)}. These filters are based on a continuous function that mimics the ideal low pass filter. The most common are Butterworth, Chebyshev and Elliptic Filters \citep{Hayes2011,Oppenheim2009}. These filters have no zeros and $N$ poles, with $N$ the filter order, equally spaced around the unit circle. 


Once these poles are computed, the continuous frequency response function can be readily obtained in its factor form and the last step consist in identifying the associated recursive formula as in exercise 4. However, note that mapping from $s$ to $z$ is usually performed using the bilinear transform\footnote{which reads: $$s=\frac{2}{\Delta t}\Biggl(\frac{1-z^{-1}}{1+z^{-1}} \Biggr)\longleftrightarrow z=\frac{1+\Delta t/2\,s}{1-\Delta t/2\,s}$$} rather than the standard mapping $z=e^{\Delta t s}$ that is used in Exercise 4, since this has the advantage of mapping $f_s/2$ to $\pi$ and thus prevent aliasing. The nonlinearity in the bilinear transform results in a \emph{wrapping} of higher frequencies so the correct cut-off frequency should first be \emph{pre-warped} to account for the distortion in the frequency calculation\footnote{The pre-warp can be achieved using $f'_c={f_s}/{\pi}\tan\bigl({\pi f_c}/{f_s}\bigr)$. Therefore, if the desired cut-off frequency is $f_c=200 Hz$ with a sampling frequency $f_s=1000Hz$, the filter should target a cut-off frequency of $f'_c=231.26 Hz$ to compensate for the warping due to the bilinear transform.}. Software packages such as \textsc{SciPy} in \textsc{Python} or \textsc{Matlab} offer the functions \textsc{butterworth} to design a Butterworth filter with given order and cut-off frequency (see \textsc{Python} script \textsc{Ex5.py}). 

The red curves in Figure \ref{fig7}a) show the amplitude response and the impulse response of a Butterworth filter of order $11$. The main advantage of these filters is their capability of well approximating the ideal filter using a limited order, which requires storing few coefficients in their recursive formulation. On the other hand, these filters tend to become \emph{unstable} as the order increases, (and the poles approach the unit circle). Moreover, their phase delay is generally not constant, and this potentially introduces phase distortion. Finally, note that since the impulse response of these filters is infinite, these cannot be implemented in the time domain via simple convolution, but via the recursive solution of the filter's LCCDE.

\textbf{2) Finite Impulse Response Filters (FIR)}. These filters are constructed in the discrete domain and have no poles (no feedback coefficients in their recursive formulation). This leads to a finite impulse response. The classic design method is the \emph{windowing} technique, which consists in multiplying the impulse response of the ideal filter in \eqref{eq42} by a window $w[k]$ which is zero outside the interval $0\leq k\leq N$, with $N$ the filter order.		

Taking $N$ as an odd number, these windows are symmetric about the midpoint, i.e. $w[k]=w[N-k]$; this results in the lag $\alpha=(N-1)/2$ in the output with respect to the input. The need for a $lag$ in the ideal filter in \eqref{eq42}-\eqref{eq43} is now clear: if $\alpha=0$, the windowed impulse response is centered in $0$, and the filter is non-causal. If the filtering is performed `offline',  it is possible to obtain a zero-phase filter by centering the windowed impulse response in $k=0$. 

Common functions are the \emph{Hanning}, \emph{Hamming} or \emph{Blackman} and \emph{Kaiser} windows. The windowing in the time domain corresponds to a convolution in the frequency domain between the ideal filter and the Fourier transform of the window function. This \emph{smooths} the transition from the band-pass to the band-stop region. Software packages as \textsc{SciPy} in \textsc{Python} or \textsc{Matlab} offer the functions \textsc{firwin} and \textsc{fir1} to design FIR filters with a given order, cut-off frequency and window function.

In addition to the linear phase response, these filters are also always stable because of the lack of poles. Moreover, the finite length of their impulse response enables their implementation via convolution. Note, however, that FIR filters require much larger order to achieve performances comparable with IIR filters. Figure \ref{fig7}b shows, in the dashed blue line, the amplitude function and the corresponding impulse response (shifted by $(N-1)/2$) of a FIR filter designed using a Hamming window of order $N=111$. The highest the filter order, the larger is the window multiplying the ideal impulse response, the more this filter approach the ideal one. On the other hand, increasing the filter order increases the sensitivity of the filter to the finite duration of the signal.

A FIR formulation make the calculation of complementary filters particularly simple thanks to the constant phase response, solely linked to the filter order. To illustrate this, consider a signal $u[k]$, and its DFT $U[n]$. Let us low-pass filter this signal to obtain $u_{\mathcal{L}}[k]$ using a FIR filter with frequency transfer function $H_{\mathcal{L}}=|H_{\mathcal{L}}|e^{-\mathrm{i}\alpha k}$. The filter operation in the frequency domain reads $U_{\mathcal{L}}=|H_{\mathcal{L}}|e^{-\mathrm{i}\alpha \theta} U$. Let $U_{\mathcal{H}}=|H_{\mathcal{H}}|e^{-\mathrm{i}\alpha \theta} U$ denote the high pass filter that gives the signal $u_{\mathcal{H}}[k]$ having complementary spectra (i.e. $|H_{\mathcal{L}}|+|H_{\mathcal{H}}|=1$) and \emph{same order} and thus \emph{same} phase delay. Using the shifting properties of the Fourier transforms, the link between the low-pass and the high-pass counter parts in the time and frequency domain sets

\begin{equation}
	\label{eq44}
	\begin{split}
		&k\,\mbox{Domain}: u_{\mathcal{H}}[k]=u[k+\alpha]-u_{\mathcal{L}}[k] \\ &\theta\,\mbox{Domain}:H_{\mathcal{H}} (\mathrm{i}\theta)U(\mathrm{i}\theta)=U(\mathrm{i}\theta) e^{\mathrm{i}\alpha \theta} -H_{\mathcal{L}}(\mathrm{i}\theta) U(\mathrm{i}\theta)
		\\&\quad \quad \quad \quad\quad\rightarrow H_{\mathcal{H}}(\mathrm{i}\theta)=e^{\mathrm{i}\alpha \theta}-H_{\mathcal{L}}(\mathrm{i}\theta)\,.
	\end{split}
\end{equation} Note that the backward shifting $u[k+\alpha]$ in the time domain cancels the phase delay produced by the low pass filter before performing the subtraction. It is easy to show that because of the linearity of the convolution, the impulse responses of complementary high-pass ($h_{\mathcal{H}}$) and low-pass ($h_{\mathcal{L}}$) filters are linked by\footnote{Note that this is not the only method to obtain an high pass filter from a low-pass filter: another approach is to reverse the frequency response $H_{\mathcal{L}}$, flipping it from left to right about the frequency $f_s/4$ for $f>0$ and from right to left about $-f_s/4$ for $f<0$ \citep{DSP}. The impulse response of the resulting high pass filter is $h_{\mathcal{H}}[n]=(-1)^n h_{\mathcal{L}}[n]$. The two methods are equivalent if the cut-off frequency separating the transition bands is $f_s/4$. This is the case encountered when performing MRA via dyadic wavelets as discussed in Chapter 5.} $\delta[n]=h_{\mathcal{L}}[n]+h_{\mathcal{H}}[n]$.

Finally, we close with the practical implementation of MRA in `off-line' conditions, for which it is possible to release the constraints of causality and use zero-phase filters. These are usually implemented by operating on the signal twice (first on $u[k]$ and then on $u[-k]$), to artificially cancel the phase delay of the operation. In \textsc{SciPy} and in \textsc{Matlab} this is performed using the function \textsc{filtfilt}.

If the phase delay is canceled, complementary filters can be computed by taking differences of the frequency transfer functions (which become real functions).
Therefore, if the first scale with band-pass $[0,f_1]$ is identified by the frequency transfer function $H_1(f)=H_{\mathcal{L}}(f,f_1)$, the second scale with band-pass $[f_1,f_2]$ is identified by a filter with transfer function $H_2(f)=H_{\mathcal{L}}(f,f_2)-H_{\mathcal{L}}(f,f_1)$. The transfer function of the general 
band-pass filter is $H_m(f)=H_{\mathcal{L}}(f,f_m)-H_{\mathcal{L}}(f,f_{m-1})$ while the last scale is identified by the high-pass filter with $H_{M}(f)=1-H_{\mathcal{L}}(f,f_{M})$.

This set of cascaded filters is known as \emph{filter bank} and is at the heart of the pyramid algorithm for computing the discrete wavelet transform \citep{Mallat2009,Strang1996}, where it is combined with sub-sampling at each scale. The general architecture of this decomposition is summarized in Figure \ref{fig7}. Observe that at the limit at which all the frequency bands become unitary, the MRA becomes a DFT.

\begin{figure}[htbp]
	\centering
	\includegraphics[keepaspectratio=true,width=1 \columnwidth]{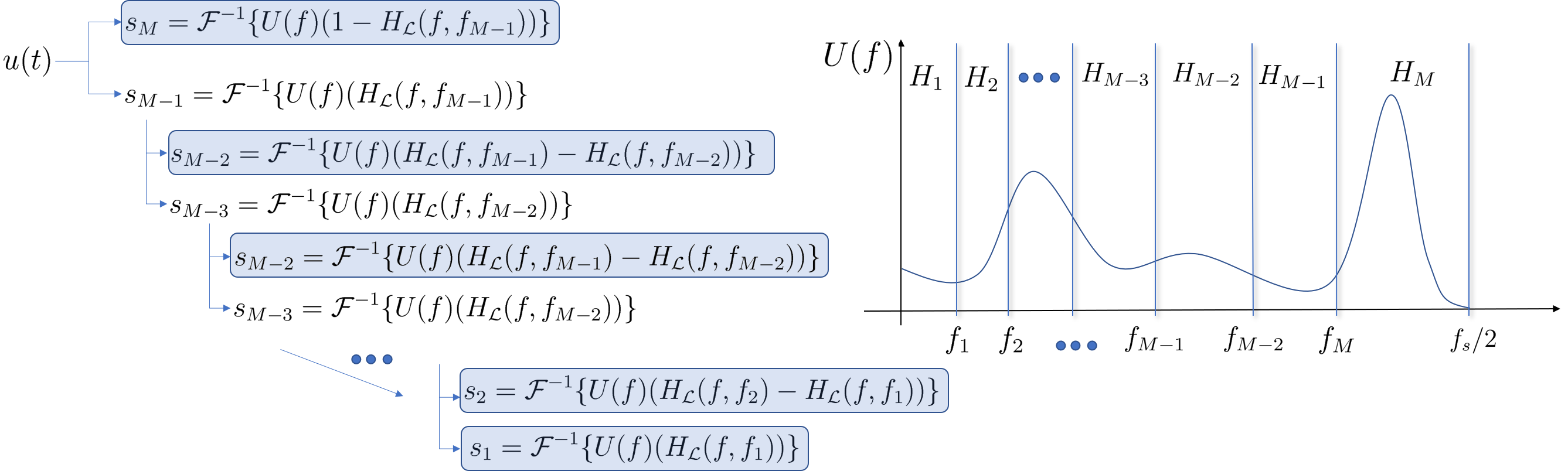}
	\caption[Examples of Discrete Signals]{Pyramid-like algorithm to compute the MRA of a signal. Each band-pass scale is computed as the difference of two low-pass filters and the terms in blue are preserved to form the summation in \eqref{eq41}. The figure on the right shows a pictorial representation of the partitioning of the signals spectra. }
	\label{fig8}
\end{figure}

\bigskip
\begin{tcolorbox}[breakable, opacityframe=.1, title=Exercise 5: Multi-Resolution Analysis]
	
	Consider the synthetic signal 
	
	\begin{equation}
		\begin{split}
			\label{eq45}
			u(t_k)=a_1 \sin(2\pi f_1 t_k) e^{\frac{-(t_k-\tau_1)^2}{b_1}}+a_2 \sin(2\pi f_2 t_k) e^{\frac{-(t_k-\tau_2)^2}{b_2}}+\\a_3\sin(2\pi f_3 t_k)+\mathcal{N}(0,0.2)
		\end{split}
	\end{equation} with $a_1=2$, $a_2=a_3=1$, $f_{1,2,3}=1, 20, 90$ (in Hz), $\tau_1=1$, $\tau_2=2.2$ (in s) and $b_1=b_2=0.05$. $\mathcal{N}(0,0.2)$ denotes Gaussian noise with zero mean and standard deviation $0.2$. The signal is sampled over $n_t=4096$ points at a sampling frequency of $f_s=1 kHz$. 
	
	1) Compute `by hand' the impulse response $h_L$ of a low pass FIR filter of order $N_O=5$, with cut off frequency $f_c=10 Hz$ using a Hamming window. Check your answers using Python's function \textsc{firwin}.
	
	2) Prepare a filter $h_L$ of order $511$ using \textsc{firwin}. Apply this to the signal using four methods. 1) direct convolution, 2) FFT-based convolution \emph{with} \textsc{Scipy}'s function \textsc{fftconvolve}, 3) FFT-based convolution using \textsc{fft} and \textsc{ifft} and 4) a recursive implementation solving the filter's LCCDE.  
	
	3) Construct a filter bank with frequency splitting vector $F_V=[10,70,100,300]$ and show the portions of the signal within the identified scales. Use filters with Hamming window and an order of $N_O=511$.

\end{tcolorbox}

\section{Application II: Time Series Analysis}\label{sec4_7}

LTI systems are the simplest model in time series analysis and forecasting. In these applications, treating signals and system as fully deterministic is too optimistic, and it is thus essential to consider stochastic signals: predictions have a certain probability range 
\citep{Guidorzi,Brockwell2010}. This section briefly reviews the main features of stochastic signals and systems in \ref{sec4_7_1}. Section \ref{sec4_7_2} reviews the basic tools for forecasting, using classic linear regression. Only the discrete domain is considered. More advanced techniques are discussed in Chapter 12.

\subsection{Stochastic LTI Systems}\label{sec4_7_1}

A stochastic signal (or the stochastic portion of a signal) is a member of an ensamble of signals characterized by a set of probability density functions.
For a comprehensive review of stochastic signals, the reader is referred to classic textbooks \citep{ljung_book,Oppenheim2015,Hsu2013}. Here, we briefly recall how LTI systems manipulates stochastic signals.

The notion and the role of the impulse response remains the same as for deterministic signals: given an input (stochastic) discrete signal, the response of the system is governed by the convolution sum in \eqref{eq21}. On the other hand, the notion of frequency spectra requires some adaptation, as stochastic signals do not generally admit a Fourier transform and focus must be placed on properties that are deterministic \emph{also} in a stochastic signal. These are the statistical properties. 

Accordingly, the time-invariance in LTI systems is extended in terms of invariance of the statistical properties. This is linked to the notion of \emph{stationarity}. Stationarity can be weak or strong. A stochastic signal $v[k]$ is stationary in a strict sense (strong stationarity) if its distributions remain invariant over time. Weak stationarity (or stationarity in a wide-sense) requires that only its \emph{time average} $\mu_v$ and \emph{autocorrelation} $r_{vv}[m]$ of a signal $v$ are time-invariant. These are defined as follows 

\begin{equation}
	\label{eq46}
	\mu_v=\mathbb{E}\{ v[k]\} \quad  \mbox{and} \quad  r_{vv}[m]= \mathbb{E}\{ v[k] v[k+m]\} \,,
\end{equation} with $\mathbb{E}$ the expectation operator. 

In the analysis of LTI system's response to stochastic signals, the link between a specific input and the corresponding output is not particularly interesting. Instead, we focus on the link between \emph{the statistical properties} of the input and the output. In particular, we consider how the properties in \eqref{eq46} are manipulated. Let $y_v[k]$ be the response of the system to the stochastic signal $u_v[k]$. The \emph{expected} (time average of the) output $\mu_{y}$ is:

\begin{equation}
	\label{eq48}
	\mu_{y}=\mathbb{E}\Biggl\{\sum^{\infty}_{l=-\infty} h[l] u_v[k-l]\Biggr\}=\sum^{\infty}_{l=-\infty} h[l] \mathbb{E}\{u_v[k-l]\} = \mu_{v}\,\sum^{\infty}_{l=-\infty} h[k]\,.
\end{equation}

This is a direct application of the homogeneity \eqref{eq16} and superposition \eqref{eq17}. We thus see that in a BIBO stable system (satisfying \eqref{eq27}) the output average is finite if the input average is finite\footnote{We also see why a high pass filter has $\sum_k h[k]=0$ while a low pass filter has $\sum_k h[k]=1$.}. The input/output relation for the the autocorrelation function has a more involved derivation, here omitted \citep{Oppenheim1996a}. Given $r_{uu}$ and $r_{yy}$ the input and the output autocorrelations, one retrieves:

\begin{equation}
	\label{eq49}
	r_{yy}[m]=\sum^{\infty}_{l=-\infty} r_{uu}[m-l] r_{hh}[l] \,\,\mbox{where}\,\, r_{hh}[k]=\sum^{\infty}_{l=-\infty} h[k] h[l+k]\,.
\end{equation}

The sequence $ r_{hh}[k]$ is the autocorrelation of the impulse response and operation on the left is a convolution. In words: \emph{the autocorrelation of the output is the convolution of the autocorrelation of the input with the autocorrelation of the impulse response}. This equation extends the convolution link in \eqref{eq21} to the autocorrelation functions. These functions admit Fourier transform, so the convolution theorem can be used to see the link in the frequency domain:

\begin{equation}
	\label{eq50}
	R_{yy}({\mathrm{j}\omega})= C_{hh} ({\mathrm{j}\omega}) R_{uu}({\mathrm{j}\omega})\,,
\end{equation} where $R_{yy}({\mathrm{j}\omega})$, $C_{hh}({\mathrm{j}\omega})$ and $R_{uu}({\mathrm{j}\omega})$ are the Fourier transform of $r_{yy}[k]$, $r_{hh}[k]$ and $r_{uu}[k]$ respectively. These are the \emph{power-spectral densities} of $y$, $h$ and $u$. Hence we see that an LTI system acts on the frequency content of the autocorrelation function of a stochastic signal.

\subsection{Time Series Forecasting via LTI Systems }\label{sec4_7_2}

Consider the explicit form of the LCCDE of a LTI in \eqref{eq37} and assume that the input signal ($u$) has both a deterministic ($u_d$) and a stochastic ($u_s$) part (i.e. $u=u_d+u_s$). Because of homogeneity, the output of the LTI system also has a deterministic ($y_d$) and a stochastic part ($y_s$). We could split these as follows:

\begin{equation}
	\label{eq51}
	y[k]=y_d[k]+y_s[k]=\sum^{N_b}_{n=0} b^*_n u_d[k-n]-\sum^{N_f}_{n=1} a^*_n y_d[k-n] + y_s\,.
\end{equation} 

Many models can be obtained depending on the assumptions on $u_s$, and hence $y_s$ \citep{Guidorzi,Nielsen2019,Brockwell2010}. For example, the stochastic part $y_s$ can be taken as white noise with zero average or as the output of a moving average filtering\footnote{A `moving average' filter is a filter with constant impulse response} of white noise. Any other filter can be used to allow controlling and/or modeling the frequency content of the stochastic contribution using \eqref{eq49} and \eqref{eq50}. 

To illustrate the main steps of time series forecasting, let us consider the simplest approach of $y_s$ being white noise. If the system is known (i.e. the coefficients $a^*_n$ and $b^*_n$ are known), the recursive equation \eqref{eq51} can be written as a matrix multiplication\footnote{This needs to be evaluated from from the first to the last entry of $\mathbf{y}^{(0)}$.}. Assume that we have collected $n_t$ samples of the input $u$.  Let $\mathbf{u}_d^{(-l)}, \mathbf{y}_d^{(-l)} \in\mathbb{R}^{n_t\times1}$ be the vectors collecting the deterministic inputs and outputs shifted backward by a lag $l$. Then, matrix form of \eqref{eq51} is

\begin{equation}
	\label{eq52}
	\begin{bmatrix}|\\\mathbf{y}^{(0)}\\|\end{bmatrix}= 
	\begin{bmatrix}
		| & | & &|& | && |\\
		\mathbf{u}^{(0)} & \mathbf{u}^{(-1)}  & \dots &\mathbf{u}^{(-N_f)}&
		\mathbf{y}^{(-1)} &\dots& \mathbf{y}^{(-N_b)} \\
		| & | & &|& |& & |
	\end{bmatrix} {\begin{bmatrix}b_0\\\vdots\\ b_{N_b}\\a_1\\ \vdots\\ a_{N_F}\end{bmatrix}}+\mathbf{y}_s.
\end{equation} 

We define $\mathbf{H}:=[\mathbf{u}^{(0)},\dots \mathbf{u}^{(-N_f)},\mathbf{y}^{(-1)}\dots \mathbf{y}^{(N_b)}]$ the Hankel matrix of the LTI system and $\mathbf{w}=[b_0,\dots b_{N_b}, a_1 \dots a_{N_f}]^T$ the vector of coefficients. The LTI system's output is

\begin{equation}
	\label{eq53}
	\mathbf{y}^{(0)}=\mathbf{H}\mathbf{w}+\mathbf{y}_s
\end{equation}

Time series forecasting via LTI systems begins with system identification. An excellent tutorial on the topic is provided by \cite{Semeraro2016}. The first goal is to identify the set of coefficients $\mathbf{w}$ from the input/output vectors. The stochastic part is considered as noise and the determinist part is our \emph{expectation}. We thus seek to solve the system $\mathbf{y}^{(0)}\approx \mathbf{H}\mathbf{w}$. Like most regression problems, this problem is ill-posed: $\mathbf{H}$ is rectangular and there is no guarantee that a unique solution exists. Like all \emph{linear} regression problems, the solution is found by minimizing a regularized cost function of the form

\begin{equation}
	\label{eq54}
	J(\mathbf{w})=||\mathbf{y}^{(0)}- \mathbf{H}\mathbf{w}||^2_2 + \alpha R(\mathbf{w})\,,
\end{equation} with $\alpha\in \mathbb{R}$ acting as a smoothing parameter and $R(\mathbf{w})$ a regularizing function\footnote{readers familiar with Lagrangian multipliers should recognize in \eqref{eq54} an augmented cost function with $\alpha$ the Lagrangian multiplier.} Classic choices are $R(\mathbf{w})=||\mathbf{w}||_2$ ($l^2$ penalty) or $R(\mathbf{w})=||\mathbf{w}||_1$ ($l^1$ penalty) or a combination of the two. 
The first is known as Tikhonov regularization, the second as LASSO regularization and the third as Elastic Net. These classic tools from machine learning \citep{Bishop2016,VladimirCherkassky2008,Murphy2012} are also employed in Chapter 12.

The reader should notice that the regression method can be generalized easily: one could replace the predictive equation \eqref{eq53} by a more complex model (e.g. an Artificial Neural Network or templates of polynomial nonlinearities as in Chapter 12), and minimize a cost function like \eqref{eq54} using an arsenal of optimization strategy. 

\bigskip

\begin{tcolorbox}[breakable, opacityframe=.1, title=Exercise 6: Time Series foreasting via LTI Systems]
	
	Generate a synthetic stochastic system with coefficients $(a_1,a_2)=(-1/2,0)$ and $(b_0,b_1,b_2)=(2,1.8,-1.2)$. Consider a deterministic input signal $u[k]=u(t_k)=\sin(8\pi t_k)+3 \exp(-(t_k-3)^2/2)$ with $t_k=k \Delta t$, $\Delta t=0.02$ and $k\in [0,1999]$. The stochastic contribution $y_s[k]$ is white noise in $[0,1]$.
	
	Use the set of 2000 points in input and output to identify the set of coefficients from the data, by solving \eqref{eq54} with a Tikhonov regularization. Consider $\alpha=1$, and recall that the minimization of $J(\mathbf{w})$ in this kind of regression (known also as Ridge regression) has a simple solution: $$\mathbf{w}=(\mathbf{H}^T\mathbf{H}+\alpha \mathbf{I})^{-1}\mathbf{H}\, \mathbf{y}^{(0)}\,,$$ with $\mathbf{I}$ the identity matrix.
	Is the expected result recovered? is the identified LTI system capable of predicting what happens if $t_k>5$? what happens for $\alpha=0$ (i.e. ordinary least square)?
\end{tcolorbox}

\section{What's next?}\label{sec4_8}

This chapter reviewed the fundamentals of signals and systems and presented Linear Time Invariant (LTI) systems in case of Single Input Single Output (SISO). We have seen that the input-output relation can be derived from knowledge of the impulse response of a system and via convolution integral. It was shown that complex exponentials are eigenfunctions of these systems and that important transforms can be derived by projecting input and output signals onto these eigenfunctions. In the eigenspace of the LTI systems, convolutions become multiplications.

Chapter 10 presents LTI systems the so-called \emph{state-space} representation, which is more common in dynamical system theory and which allows for straightforward generalization to MIMO systems. 
Chapter 11 reviews the analysis of nonlinear systems, while Chapter 12 describes the  system identification more broadly, considered also nonlinear systems.

This chapter reviewed the link between the continuous and the discrete world and the impact of the discretization on the eigenfunctions of an LTI system. Special values of the complex frequencies, called poles, yields infinite response of a system and are linked to the notion of stability, reviewed in Chapter 13. The reader should recognized that the identification of these poles from large datasets is the essence of the Dynamic Mode Decomposition described in Chapter 7. 

Finally, this chapter also introduced the fundamentals of MRA, which well complements wavelet theory in Chapter 5 and Chapter 8 on the multiscale Proper Orthogonal Decomposition.


\bibliographystyle{apalike}
\bibliography{Chapter_4_DDFM_Book_Mendez}
\end{document}